\documentclass[notitlepage,11pt]{article}
\usepackage{amssymb,amsfonts,amsmath,geometry,esint,comment,upgreek,cases,
graphicx,ulem}
\usepackage{scalerel,stackengine,mathtools,stmaryrd}

\usepackage{color}

\catcode`\@=11
\@addtoreset{equation}{section}

\catcode`\@=12

\allowdisplaybreaks

\def\eps{\varepsilon}

\def\R{\mathbb{R}}

\def\f{\varphi}

\def\S{\mathbb S}

\def\ird{\int\limits_{\R^d}}

\def\isd{\int\limits_{\S^{d-1}}}

\def\C{\mathcal C}

\def\div{{\rm div}}

\def\rmH{{\rm H}}

\def\proof{\noindent{\textbf{Proof. }}}
\def\QED{\hfill {$\square$}\goodbreak \medskip}

\newtheorem{Theorem}{Theorem} 
\newtheorem{Lemma}{Lemma}

\newtheorem{Corollary}{Corollary} 
\newtheorem{Remark}[Theorem]{Remark}

\begin{document}

\title{
Hardy type inequalities with mixed weights in cones}

\author{{Gabriele Cora\footnote{Dipartimento di Matematica "G. Peano", Universit\`a di Torino, Italy. Email: {gabriele.cora@unito.it}, orcid.org/0000-0002-0090-5470.}}, 
Roberta Musina\footnote{Dipartimento di Scienze Matematiche, Informatiche e Fisiche, Universit\`a di Udine, Italy.
Email: {roberta.musina@uniud.it}, orcid.org/0000-0003-4835-8004.},
Alexander I. Nazarov\footnote{
St. Petersburg Department of Steklov Institute and St{.} Petersburg State University, St{.} Petersburg, Russia. E-mail: al.il.nazarov@gmail.com, 
orcid.org/0000-0001-9174-7000.}
}

\date{}

\maketitle

\centerline{\it Dedicated to Susanna Terracini}

\bigskip

\noindent
{\small {\bf Abstract.} We study Hardy type inequalities involving mixed cylindrical and spherical weights, for functions supported in cones. These inequalities are related to some singular or degenerate 
differential operators.}



\section{Introduction}

We deal with the best constant $m_{p,a,{b}}(\mathcal C)$ in inequalities of the type
\begin{equation}
\label{eq:H_best}
m_{p,a,{b}}(\mathcal C)\int\limits_{\mathcal C}|y|^{a}|z|^{-{{b}}-p}|u|^p~\!dz\le \int\limits_{\mathcal C}|y|^{a}|z|^{-{{b}}}|\nabla u|^p~\!dz~,
\quad u\in  C^\infty_c({\mathcal C})~\!,
\end{equation}
where ${\mathcal C}\subseteq \R^d$ is a cone, that is a dilation-invariant open set, $p>1$,  $a, {b}\in\R$
and $z=(x,y)$ is the variable in $\R^d\equiv\R^{d-k}\times\R^k$.
Inequality (\ref{eq:H_best}) includes the well known case of purely spherical weights; 
for this reason we assume that $1\le k<d$ and $a\neq 0$.

Our starting motivation arose from the growing interest, inspired by \cite{CS}, in differential operators of the form
\begin{equation}
\label{eq:operator}
\mathcal Lu=-\div(|y|^aA(z)\nabla u)~\!.
\end{equation}
Indeed, in case $p=2$, the right hand side in (\ref{eq:H_best}) is the quadratic form associated to the 
differential operator $\mathcal L$, where $A(z)=|z|^{-b}~\!\mathcal{I}_d$. 
Starting with the seminal paper \cite{KFS}, large efforts have been spent to investigate degenerate/singular 
operators including (\ref{eq:operator}) (see for instance   \cite{DFM, DM} and references therein). 
We refer also to the 
papers \cite{STV1, STV2, TTV} by Susanna Terracini and collaborators. With respect to \cite{STV1}
(even solutions), and \cite{STV2}  (odd solutions), the relevant cones are $\R^d$ and $\R^d_+=\R^{d-1}\times(0,\infty)$,
respectively.

In dealing with (\ref{eq:H_best}), one is forced to assume that the weights involved are locally integrable on the cone $\mathcal C$.
This leads us to distinguish different situations, depending on the position of $\mathcal C$ with respect to the singular set 
$$
\Sigma_0:=\{y=0\}. $$
More precisely, the  cone $\mathcal C=\R^d$ needs $k+a>0$ and 
$d+a>p+b$; if $\mathcal \C\subseteq \R^d\setminus\{0\}$ then we have to require that 
\begin{equation}
\label{eq:C_assumption}
\mathcal C\subseteq\R^d\setminus\Sigma_0~,\quad\text{or}\quad k+a>0.
\end{equation}

In any case,
a special role is played by the real constant
\begin{equation}
\label{eq:def_H}
\rmH_{p,a,{b}}=\frac{d+a}{p}-\frac{p+{b}}{p}~\!.
\end{equation}
We start with a simple result, that deals with the two largest cones (here, $\partial_r $ stands for the radial derivative).

\begin{Theorem}
\label{T:H}
Let $k+a>0$. Then the inequality 
\begin{equation}
\label{eq:H_best_rad}
|\rmH_{p,a,{b}}|^p\ird|y|^{a}|z|^{-{{b}}-p}|u|^p~\!dz\le\ird|y|^{a}|z|^{-{{b}}}|\partial_r u|^p~\!dz
\end{equation}
holds for any $u\in  C^\infty_c(\R^{d}\setminus\{0\})$, with a sharp constant in the left hand side. 

In particular,
$m_{p,a,{b}}(\R^d\setminus\{0\})=|\rmH_{p,a,{b}}|^p $.

 If $d+a>p+{b}$, then (\ref{eq:H_best_rad}) holds for any $u\in  C^\infty_c(\R^{d})$, 
and $m_{p,a,{b}}(\R^d)=\rmH^p_{p,a,{b}} $. 
\end{Theorem}

Theorem \ref{T:H} implies that the Hardy inequality (with mixed weights) in $\R^d\setminus\{0\}$  holds with a positive constant
if and only if $d+a\neq p  +{b}$. From the proof,
see Section \ref{S:proofs}, it is evident that the corresponding best constant is not achieved on any reasonable function space.

\medskip

We point out a  remarkable case. It is well known that
\begin{equation}
\label{eq:cylindrical}
\Big(\frac{k+a-p}{p}\Big)^p\int\limits_{\R^d}|y|^{a-p}|u|^p~\!dz\le \int\limits_{\R^d}|y|^{a}|\nabla u|^p~\!dz
\quad \text{for any $u\in C^\infty_c(\R^d)$,}
\end{equation}
provided that $a> p-k$, which is needed for the local integrability of the weight in the left-hand side. Notice that in the threshold case
$a=p-k$, the constant in (\ref{eq:cylindrical}) vanishes. In contrast, a Hardy inequality involving {\it mixed weights} holds, as pointed out in the next statement.

\begin{Corollary}
Let $d>k$. Then
$$
\Big(\frac{d-k}{p}\Big)^p\int\limits_{\R^d}|y|^{p-k}|z|^{-p}|u|^p~\!dz\le \int\limits_{\R^d}|y|^{p-k}|\nabla u|^p~\!dz
\quad \text{for any $u\in C^\infty_c(\R^d)$.}
$$
The constant in the left hand side can not be improved. 
\end{Corollary}

From now on we deal with general cones $\mathcal C\subset\R^d$. We start with the 
{\it superdegenerate case} (cf. \cite{TTV}).

\begin{Theorem}
\label{T:old}
Let $\mathcal C\subseteq \R^d\setminus \{0\}$. If $k+a\ge p$, then $m_{p,a,b}(\mathcal C\setminus\Sigma_0)= m_{p,a,b}(\mathcal C)$. In particular,
$$m_{p,a,{b}}(\R^d\setminus\Sigma_0)=|\rmH_{p,a,{b}}|^p.$$
\end{Theorem}

The next result is already known in case of purely spherical weights \cite[Section 6]{Nsur}; see also \cite{Scone} for related issues in case $p=2$.
Here we denote by $\Pi$ the restriction of the orthogonal projection $\Pi:\R^{d-k}\times\R^k\to \{0\}\times\R^k$
to the unit sphere $\S^{d-1}$.

\begin{Theorem}
\label{T:cone_sphere1}
Let $\mathcal C$ satisfy (\ref{eq:C_assumption}) and put $\omega_\mathcal C=\S^{d-1}\cap \mathcal C$. Then $m_{p,a,{b}}(\mathcal C)=\mathcal M_{p,a,{b}}(\omega_\mathcal C)$, where
$$
\mathcal M_{p,a,{b}}(\omega_\mathcal C)=\inf_{\f\in  C^\infty_c(\omega_\mathcal C)\atop \f\neq 0}
\frac{\displaystyle\int_{\omega_\mathcal C}|\Pi \sigma|^a(|\nabla_{\!\sigma} \f|^2+{ \rmH^2_{p,a,{b}}}|\f|^2)^\frac{p}{2}d\sigma}
{\displaystyle\int_{\omega_\mathcal C}|\Pi \sigma|^a|\f|^pd\sigma}~\!.
$$
\end{Theorem}

Under the assumption (\ref{eq:C_assumption}), we can define the following weighted Sobolev spaces:
\begin{itemize}
\item[$S1)$] ${\mathcal D}^{1,p}_0(\mathcal C;|y|^a|z|^{-b}dz)$ is the completion of  $ C^\infty_c(\mathcal C)$ with respect to the norm
$$
\|u\|^p=\int\limits_\mathcal C|y|^a|z|^{-b}\big(|\nabla u|^p+|z|^{-p}|u|^p\big)~\!dz~\!;
$$
\item[$S2)$] $W^{1,p}_0(\omega_\mathcal C;|\Pi\sigma|^ad\sigma)$ is the completion of  $C^\infty_c(\omega_\mathcal C)$ with respect to the norm
$$
|\!|\!|\f|\!|\!|^p =\int\limits_{\omega_\mathcal C}|\Pi\sigma|^a\big(|\nabla_{\!\sigma} \f|^p+|\f|^p\big)~\!d\sigma.
$$
\end{itemize}

The best constant $\mathcal M_{p,a,b}({\omega_\mathcal C})$ is attained in  
$W^{1,p}_0({\omega_\mathcal C};|\Pi\sigma|^ad\sigma)$, due to the compactness of the
embedding $W^{1,p}_0(\omega_\mathcal C;|\Pi\sigma|^ad\sigma)\hookrightarrow L^p(\omega_\mathcal C;|\Pi\sigma|^ad\sigma)$ given by Corollary \ref{C:HS_sphere} in Section \ref{S:preliminary}. 
 In contrast, the next result holds.

\begin{Theorem}
\label{T:not_achieved}
Let $\mathcal C$ satisfy (\ref{eq:C_assumption}).
The infimum
$m_{p,a,b}(\mathcal C)$ is  not achieved on  ${\mathcal D}^{1,p}_0(\mathcal C;|y|^a|z|^{-b}dz)$.
\end{Theorem}

Theorem \ref{T:cone_sphere1} shows, in particular, that $m_{p,a,{b}}(\mathcal C)\ge |\rmH_{p,a,{b}}|^p$ (notice that for $k+a>0$, this evidently follows from 
Theorem \ref{T:H}) and that equality might occur, compare with Theorem \ref{T:old}.
It is natural to look for conditions that guarantee the validity of the strict inequality.

Since
$\mathcal M_{p,a,{b}}(\omega_\mathcal C)$ is achieved, then
$m_{p,a,{b}}(\mathcal C)> |\rmH_{p,a,{b}}|^p$, unless nontrivial constant functions are in 
$W^{1,p}_0(\omega_\mathcal C;|\Pi\sigma|^ad\sigma)$. Consider for instance the following situations:
\begin{itemize}
\item[$i)$] $\Sigma_0$ does not intersect $\overline{\omega}_\mathcal C$.
Then $|\Pi\sigma|^a$ is bounded and bounded away from $0$, so we are in fact in the case $a=0$; 
\item[$ii)$] $k+a<kp$. Then $W^{1,p}_0(\omega_\mathcal C;|\Pi\sigma|^ad\sigma)$ is  embedded into the standard Sobolev space $W^{1,1}_0(\omega_\mathcal C)$.
\end{itemize}
In both cases, we infer that $m_{p,a,{b}}(\mathcal C)> |\rmH_{p,a,{b}}|^p$ if, for instance,
$\overline{\omega}_\mathcal C\subsetneq \S^{d-1}$. 
More refined sufficient conditions to get $1\notin W^{1,p}_0(\omega_\mathcal C;|\Pi\sigma|^ad\sigma)$
can be obtained by using the results in  \cite[Subsection 13.2]{Ma}.

The case $k+a\ge kp$ is included in the next statement.

\begin{Theorem}
\label{T:strict1}
Let $\mathcal C$ satisfy (\ref{eq:C_assumption}).
If $\partial\mathcal C\setminus\Sigma_0\neq \emptyset$, then
$$m_{p,a,{b}}(\mathcal C)> |\rmH_{p,a,{b}}|^p.$$
\end{Theorem}

It not easy to calculate $m_{p,a,{b}}(\omega_\mathcal C)$ for general cones and exponents. 
If $a=0$ then Theorem \ref{T:cone_sphere1} gives 
$$
m_{2,0,{b}}(\mathcal C)=
\inf_{\f\in  H^1_0(\omega_\mathcal C)\atop \f\neq 0}
\frac{\displaystyle\int_{\omega_\mathcal C}(|\nabla_{\!\sigma} \f|^2+{ \rmH^2_{2,0,{b}}}|\f|^2)d\sigma}
{\displaystyle\int_{\omega_\mathcal C}|\f|^2d\sigma}=\lambda_1(\omega_\mathcal C)+|\rmH_{2,0,{b}}|^2,
$$
where $\lambda_1(\omega_\mathcal C)$ is the first eigenvalue of the Laplace-Beltrami operator on $\omega_\mathcal C$
with Dirichlet boundary conditions.

For  the cone $\R^{d}\setminus\Sigma_0$ we have the following result (here $s^+=\max\{s,0\}$).
\begin{Theorem}
\label{T:strict3}
Let $p=2$. Then $m_{2,a,{b}}(\R^{d}\setminus\Sigma_0)>|\rmH_{2,a,{b}}|^2$ if and only if $k+a<2$. More precisely,
$$
m_{2,a,{b}}(\R^{d}\setminus\Sigma_0)=(d-k)(2-(k+a))^++ |\rmH_{2,a,{b}}|^2.
$$
\end{Theorem}

In case $k=1$, the singular set $\Sigma_0=\{y=0\}$ is a hyperplane which disconnects $\R^d$. Thus 
from Theorems \ref{T:old}, \ref{T:strict3} we immediately obtain the next statement.

\begin{Corollary}
\label{C:cor}
\begin{itemize}
\item[$i)$]
If $a\ge p-1$, then $m_{p,a,b}(\R^d_+)= |{\rmH_{p,a,b}}|^p$. 
\item[$ii)$]
If $p=2$, then $m_{2,a,b}(\R^d_+)= (d-1)(1-a)^++ |\rmH_{2,a,{b}}|^2$ for any $a\in\R$.
\end{itemize}
\end{Corollary}

In the Caffarelli-Silvestre \cite{CS}  setting   we have $d=n+1$, $a=1-2s\in (-1,1)$ and $b=0$. 
In case $n>2s$ (which is a restriction only if $n=1$), the equality
$$
\inf_{u\in C^\infty_c(\R^{n+1})\atop \f\neq 0}
\frac{\displaystyle\int_{\R^{n+1}}|y|^{1-2s}|\nabla u|^2dz}
{\displaystyle\int_{\R^{n+1}}{|y|^{1-2s}}{|z|^{-2}}| u|^2dz}=\Big(\frac{n-2s}{2}\Big)^2
$$
(see for instance \cite[formula (11)]{MN_AiHP}), follows via Corollary \ref{C:cor} as well. On the  half space $\R^{n+1}_+$ we have  
$$
\inf_{u\in C^\infty_c(\R^{n+1}_+)\atop \f\neq 0}
\frac{\displaystyle\int_{\R^{n+1}_+}|y|^{1-2s}|\nabla u|^2dz}
{\displaystyle\int_{\R^{n+1}_+}{|y|^{1-2s}}{|z|^{-2}}| u|^2dz}=\Big(\frac{n+2s}{2}\Big)^2~\!.
$$
Notice that no restriction on $s$ in case $n=1$ is needed.

\medskip
The paper is organized as follows. The preliminary Section \ref{S:preliminary} is mostly 
devoted to the properties of the weighted Sobolev space $W^{1,p}_0(\omega;|\Pi\sigma|^ad\sigma)$.
The proofs of the theorems  stated above are collected in Section \ref{S:proofs}.

\bigskip

{\small
\noindent
{\bf Notation.}
We denote by $\Pi$ both the orthogonal projection $\Pi:\R^d\to \{0\}\times\R^k$,
and its restriction to the unit sphere $\S^{d-1}$. 
So, if $r>0, \sigma\in\S^{d-1}$ are the spherical coordinates of $z=(x,y)\in \R^{d-k}\times\R^k$, then
$|y|=r|\Pi\sigma|$.

For $R>0$ we denote by $B_R$ the ball of radius $R$ about the origin.

Through the paper, any positive constant whose value is not important is denoted by $c$. It may take different values at different places. To indicate that a constant depends on some parameters 
we list them in parentheses. 
}

\section{Preliminaries}
\label{S:preliminary}

We start by pointing out  the local integrability 
properties of the weights involved.

\begin{Lemma}
\label{L:L1loc}
Let ${{a}},\beta\in\R$. Then 
\begin{itemize}
\item[$i)$] 
$|y|^{{a}}|z|^{-\beta}\in L^1_{\rm loc}(\R^d\setminus\{0\})$ if and only if 
$k+{{a}}>0$;
\item[$ii)$] $|\Pi \sigma|^a\in L^1(\S^{d-1})$ if and only if $k+a>0$;
\item[$iii)$] $|y|^{{a}}|z|^{-\beta}\in L^1_{\rm loc}(\R^d)$ if and only if 
$k+{{a}}>0$ and $d+{{a}}>\beta$.
\end{itemize}
\end{Lemma}

\proof
Evidently, $|y|^{{a}}|z|^{-\beta}\in L^1_{\rm loc}(\R^d\setminus\{0\})$ if and only if $|y|^{a}\in L^1_{\rm loc}(\R^k)$, that is, $k+a>0$. 
Since 
$$
\int\limits_{B_R}|y|^{{a}}|z|^{-\beta}~\!dz=\int\limits_0^R r^{d+a-\beta-1}dr\int\limits_{\S^{d-1}}|\Pi\sigma|^a~\!d\sigma
$$
for any $R>0$, then $ii)$ and $iii)$ readily follow.
\QED

The main results in this section deal with the weighted Sobolev spaces $W^{1,p}_0(\omega;|\Pi\sigma|^ad\sigma)$ for $\omega\subseteq\S^{d-1}$ open, under the  assumption
\begin{equation}
\label{eq:omega_ass}
 \omega\cap\Sigma_0=\emptyset\qquad \text{if $k+a\le 0$.}
\end{equation}
Since $|\Pi\sigma|^a\in L^1_{\rm loc}(\omega)$, the space $W^{1,p}_0(\omega;|\Pi\sigma|^ad\sigma)$ is plainly defined,
accordingly with $S2)$ in the Introduction.

\begin{Lemma} 
\label{L:continuous}
Let $\omega$ satisfy (\ref{eq:omega_ass}). If 
$\S^{d-1}\cap \Sigma_0\not\subset \overline\omega$,
then
there exist $t>0$ depending only on $a,k$ and $p$, and $c(\omega)>0$ such that 
\begin{equation}
\label{eq:ball}
\displaystyle
\int\limits_{\omega\cap\{|\Pi\sigma|<\eps\}}\!\!\!\!|\Pi\sigma|^a|\f|^p~\!d\sigma\le 
c(\omega)\eps^t\int\limits_{\omega}|\Pi\sigma|^a|\nabla_{\!\sigma} \f|^p~\!d\sigma
\end{equation}
for any $\f \in  W^{1,p}_0(\omega;|\Pi\sigma|^ad\sigma)$ and any  $\eps>0$.
 \end{Lemma}

\proof
We can assume that the south  pole ${\rm e}=(-1,0\dots,0)$ does not belong to $\overline\omega$.

Let ${\rm P}:\R^{d-1}\equiv  \R^{d-1-k}\times\R^k\to \S^{d-1}\setminus\{{\rm e}\}\subset \R^{d}$ be the inverse of the stereographic projection from ${\rm e}$.
More explicitly, 
$$
{\rm P}(\xi,y)=(\mu-1,\mu \xi;\mu y)\in \S^{d-1}\subset (\R\times\R^{d-k-1})\times\R^k~,\quad \mu(\xi,y)=\frac{2}{1+|\xi|^2+|y|^2}
$$
(the variable $\xi$ has to be omitted if $k=d-1$). 

Let $\Omega:={\rm P}^{-1}(\omega)$. Then $\Omega\subset\R^{d-1}$ is open and bounded. Moreover,
if $k+a\le 0$ then $\Omega\cap \Sigma_0= \emptyset$ because  $\omega\cap \Sigma_0= \emptyset$ by (\ref{eq:omega_ass}). 

Next, for $\f\in C^\infty_c(\omega)$ we put 
$\tilde{\f}:=\f\circ {\rm P}\in  C^\infty_c(\Omega)$. Since $\mu=\mu(\xi,y)$ is bounded and bounded away from $0$ on $\Omega$, 
we have that there exist $c>1, R>0$ depending only on $\omega$ such that
\begin{equation}
\label{eq:mu}
\begin{aligned}
\int\limits_{\omega}|\Pi \sigma|^a|\nabla_{\!\sigma} \f|^p~\!d\sigma&=
\int\limits_{\Omega}\mu^{d+a-p-1}|y|^a|\nabla {\tilde{\f}}|^p~\!d\xi dy\ge
c^{-1}\int\limits_{\Omega}|y|^a|\nabla {\tilde{\f}}|^p~\!d\xi dy
\\
\int\limits_{\omega\cap\{|\Pi\sigma|<\eps\}}|\Pi \sigma|^a|\f|^pd\sigma&\le
\int\limits_{\Omega\cap\{|y|<R\eps\}}\mu^{d+a-1}|y|^a|{\tilde{\f}}|^pd\xi dy
\le c \int\limits_{\Omega\cap\{|y|<R\eps\}}|y|^a|{\tilde{\f}}|^pd\xi dy.
\end{aligned}
\end{equation}

The conclusion of the proof follows via the Hardy-Maz'ya inequalities in \cite[Chapter 2]{Ma}. 

Let $k+a>0$. We fix any $0<t<\min\{k+a,p\}$ and use 
$$
\int\limits_{\R^{d-1}}|y|^{a-t}|\tilde\f|^p~\!d\xi dy\le \big(\tfrac{p}{k+a-t}\big)^{p} \int\limits_{\R^{d-1}}|y|^{a+p-t}|\nabla \tilde\f|^p~\!d\xi dy~\!,
$$
to estimate
$$
\int\limits_{\Omega\cap\{|y|<R\eps\}}\!\!\!\!|y|^a|\tilde\f|^p~\!d\xi dy\le c\eps^t \int\limits_{\R^{d-1}}|y|^{a-t}|\tilde\f|^p~\!d\xi dy\le c \eps^t \int\limits_{\R^{d-1}}|y|^{a+p-t}|\nabla \tilde\f|^p~\!d\xi dy
\le c \eps^t \int\limits_{\Omega}|y|^{a}|\nabla \tilde\f|^p~\!d\xi dy~\!.
$$
Thus (\ref{eq:ball}) is proved in this case, thanks to (\ref{eq:mu}).

If $k+a\le 0$ we use
$$
\int\limits_{\R^{d-1}}|y|^{a-p}|\tilde\f|^p~\!d\xi dy\le \big(\tfrac{p}{p-(k+a)}\big)^{p} \int\limits_{\R^{d-1}}|y|^a|\nabla \tilde\f|^p~\!d\xi dy
$$
to get
$$
\int\limits_{\Omega\cap\{|y|<R\eps\}}\!\!\!\!|y|^a|\tilde\f|^p~\!d\xi dy\le c\eps^p \int\limits_{\R^{d-1}}|y|^{a-p}|\tilde\f|^p~\!d\xi dy\le c \eps^p \int\limits_{\R^{d-1}}|y|^a|\nabla \tilde\f|^p~\!d\xi dy,
$$
and (\ref{eq:ball}) again follows.
The proof is complete.
\QED

\begin{Lemma}
\label{L:compact}
Let $\omega$ satisfy (\ref{eq:omega_ass}). If 
$\S^{d-1}\cap \Sigma_0\not\subset \overline\omega$, then 
\begin{equation}
\label{eq:positive}
\int\limits_{\omega}|\Pi\sigma|^a|\f|^p~\!d\sigma\le c\int\limits_{\omega}|\Pi\sigma|^a|\nabla_{\!\sigma} \f|^p~\!d\sigma
\quad\text{for any}\quad \f\in W^{1,p}_0(\omega;|\Pi\sigma|^ad\sigma)
\end{equation}
and 
the embedding $W^{1,p}_0(\omega;|\Pi\sigma|^ad\sigma)\hookrightarrow L^p(\omega;|\Pi\sigma|^ad\sigma)$ is compact.
\end{Lemma}

\proof
Trivially, inequality (\ref{eq:positive}) follows from (\ref{eq:ball}), by choosing $\eps=1$.

The embedding operator $W^{1,p}_0(\omega;|\Pi\sigma|^ad\sigma)\hookrightarrow 
L^p(\omega \cap \{|\Pi\sigma|>\eps\};|\Pi\sigma|^ad\sigma)$ is compact by the Rellich theorem. 
Lemma \ref{L:continuous} shows that the operator 
$W^{1,p}_0(\omega;|\Pi\sigma|^ad\sigma)\hookrightarrow L^p(\omega;|\Pi\sigma|^ad\sigma)$ can be approximated in norm by compact operators. Thus it is compact itself. \QED

\begin{Corollary}
\label{C:HS_sphere} 
\begin{itemize}
\item[$i)$] The embedding $W^{1,p}_0(\S^{d-1}\setminus\Sigma_0;|\Pi\sigma|^ad\sigma)\hookrightarrow 
L^p(\S^{d-1};|\Pi\sigma|^ad\sigma)$ is compact;
\item[$ii)$] If $k+a>0$, the embedding $W^{1,p}(\S^{d-1};|\Pi\sigma|^ad\sigma)\hookrightarrow L^p(\S^{d-1};|\Pi\sigma|^ad\sigma)
$ is compact.
\end{itemize}
\end{Corollary}

\proof 
Take a point ${\rm e}\in\S^{d-1}\cap \Sigma_0$ and a cut-off function $\eta\in C^\infty_c(\S^{d-1}\setminus\{{\rm e}\})$
such that $\eta\equiv 1$ in a neighborhood of ${\rm -e}$. The operators
$$
\f\mapsto \eta\f~,\quad \f\mapsto (1-\eta)\f~,\qquad  W^{1,p}_0(\S^{d-1}\setminus\Sigma_0;|\Pi\sigma|^ad\sigma)\to 
L^p(\S^{d-1};|\Pi\sigma|^ad\sigma)
$$
are compact  by Lemma \ref{L:compact}, which proves $i)$. For  $ii)$ repeat the same argument.
\QED

\section{Proofs}
\label{S:proofs}

\paragraph{\bf Proof of Theorem \ref{T:H}.}
Notice that $|\Pi\sigma|^a\in L^1(\S^{d-1})$, see Lemma \ref{L:L1loc}.
Let $u\in  C^\infty_c(\R^d)$. If $d+a\le p+{b}$ assume in addition that 
$u\in  C^\infty_c(\R^d\setminus \{0\})$. 

We use the classical Hardy inequality for functions of one variable, 
which holds with a sharp and not achieved constant, see \cite[Theorem 330]{HLP}, to estimate
$$
\int\limits_0^\infty r^{d+a-{b}-1}|(\partial_r u)(r\sigma)|^p~\!dr\ge 
|\rmH_{p,a,{b}}|^p\int\limits_0^\infty r^{d+a-{b}-p-1}| u(r\sigma)|^p~\!dr
$$
for any $\sigma\in\S^{d-1}$. It follows that
$$
\begin{aligned}
\ird|y|^{a}|z|^{-{{b}}}|\partial_r u|^p~\!dz&=\isd|\Pi\sigma|^a~\!d\sigma\int\limits_0^\infty r^{d+a-{b}-1}|(\partial_r u)(r\sigma)|^p~\!dr
\\&\ge  |\rmH_{p,a,{b}}|^p\!\!\!\isd|\Pi\sigma|^a~\!d\sigma\int\limits_0^\infty r^{d+a-{b}-p-1}| u(r\sigma)|^p~\!dr=|\rmH_{p,a,b}|^p\!\!\! \ird|y|^{a}|z|^{-b-p}|u|^p~\!dz,
\end{aligned}
$$
which concludes the proof. \QED

\paragraph{\bf Proof of Theorem \ref{T:old}.}
Clearly $m_{p,a,{b}}(\mathcal C\setminus\Sigma_0)\ge  m_{p,a,{b}}(\mathcal C)$. To prove the opposite inequality
fix any nontrivial function $u\in  C^\infty_c(\mathcal C)$. Our aim is to suitably approximate $u$  by a sequence of
functions in $ C^\infty_c(\C\setminus\Sigma_0)$. We take a  function $\eta\in  C^\infty_c(\R)$ such that $\eta\equiv 1$ on $[0,1]$, $\eta\equiv 0$ on $[2,\infty)$ and $\|\eta\|_\infty\le 1$. For any integer $h\ge 1$ we
put
$$
u_h(x,y)=\eta\Big(\frac{-\log |y|}{h}\Big)u(x,y).
$$
Then $u_h\in  C^\infty(\R^k\setminus\Sigma_0)$, $u_h\to u$ pointwise and
$$
u_h(x,y)=0~~\text{if $|y|\le e^{-2h}$}~,\qquad u_h(x,y)=u(x,y)~~\text{if $|y|\ge e^{-h}$.}
$$
We claim that
\begin{equation}
 \label{eq:grad}
 \begin{aligned}
 \int\limits_{\mathcal C}|y|^{a}|z|^{-{{b}}-p}|u_h|^p~\!dz&=\int\limits_{\mathcal C}|y|^{a}|z|^{-{{b}}-p}|u|^p~\!dz+o_h(1)\\
 \int\limits_{\mathcal C}|y|^{a}|z|^{-{{b}}}|\nabla u_h|^p~\!dz&=\int\limits_{\mathcal C}|y|^{a}|z|^{-{{b}}}|\nabla u|^p~\!dz+o_h(1)
 \end{aligned}
\end{equation}
as $h\to\infty$. The first limit in (\ref{eq:grad}) plainly follows by Lebesgue's theorem. To prove the second one,
it suffices to show that
$$
I_h:= \int\limits_{S_h}|y|^{a}|z|^{-{{b}}}|\nabla u_h|^p~\!dz=o_h(1)~,\qquad S_h:=\R^{d-k}\times\{e^{-2h}<|y|<e^{-h}\}.
$$
Let $\delta\in(0,1)$ be such that $\text{supp}(u)\subset \{\delta<|z|<\delta^{-1}\}$. Since
 $$
|\nabla u_h(x,y)|\le \|\nabla \eta\|_\infty \frac{|u(x,y)|}{h|y|}+|\nabla u(x,y)| \quad\text{on $S_h$,}
$$
we can estimate
$$
\begin{aligned}
I_h&\le \frac{c}{h^p}\int\limits_{S_h}|y|^{a-p} |z|^{-b}|u|^p~\!dz+o_h(1)
\le \frac{c(\delta,u)}{h^p} \int\limits_{e^{-2h}}^{e^{-h}} r^{k+a-p-1}dr+o_h(1)=o_h(1).
\end{aligned}
$$

In conclusion, we have
$$
m_{p,a,{b}}(\mathcal C\setminus\Sigma_0)\le \frac
{\displaystyle{\int\limits_{\mathcal C}|y|^{a}|z|^{-{{b}}}|\nabla u_h|^p~\!dz}}
{\displaystyle{\int\limits_{\mathcal C}|y|^{a}|z|^{-{{b}}-p}|u_h|^p~\!dz} }
=
\frac
{\displaystyle{\int\limits_{\mathcal C}|y|^{a}|z|^{-{{b}}}|\nabla u|^p~\!dz}}
{\displaystyle{\int\limits_{\mathcal C}|y|^{a}|z|^{-{{b}}-p}|u|^p~\!dz} }+o_h(1).
$$
Since $u\in  C^\infty_c(\mathcal C)$ was arbitrarily chosen, the  inequality $m_{p,a,{b}}(\mathcal C\setminus\Sigma_0)
\le m_{p,a,{b}}(\mathcal C)$ follows. The proof of Theorem \ref{T:old} is complete.
\QED

\begin{Remark}
The same argument shows that 
$ W^{1,p}_0(\omega\setminus\Sigma_0;|\Pi\sigma|^ad\sigma)=W^{1,p}_0(\omega;|\Pi\sigma|^ad\sigma)$
for any open $\omega\subseteq\S^{d-1}$, provided that $k+a\ge p$.
\end{Remark}

\paragraph{Proof of Theorem \ref{T:cone_sphere1}.}
We have to compare the infima
$$
m(\mathcal C)=\inf_{u\in   C^\infty_c(\mathcal C)\atop u\neq 0}
\frac{\displaystyle\int_{\mathcal C}|y|^a|z|^{-b}|\nabla u|^p~\!dz}{\displaystyle\int_{\mathcal C}|y|^a|z|^{-b-p}| u|^p~\!dz}\
~,\quad
\mathcal M({\omega_\mathcal C})=\inf_{\f\in  W^{1,p}_0({\omega_\mathcal C};|\Pi\sigma|^ad\sigma)\atop \f\neq 0}
\frac{\displaystyle\int_{\omega_\mathcal C}|\Pi\sigma|^a(|\nabla_{\!\sigma}\f|^2+\rmH^2|\f|^2)^\frac{p}{2}~\!d\sigma}{\displaystyle\int_{\omega_\mathcal C} |\Pi\sigma|^a|\f|^p~\!d\sigma},
$$
where $\rmH:=\rmH_{p,a,b}$ is given by (\ref{eq:def_H}) (in this proof we omit  the indexes $p,a,b$).

\medskip

Thanks to Theorem \ref{T:old}, we can assume that the following stronger hypothesis hold:
\begin{equation}
\label{eq:stronger}
\mathcal C\subset\R^d\setminus\Sigma_0~,\quad\text{or}\quad 0<k+a<p.
\end{equation}
In the first case, the weight $|y|^a|z|^{-b}$ is bounded and bounded away from zero on 
any compact set in $\mathcal C$. If $0<k+a<p$ then $|y|^a$ belongs to the Muckenhoupt class $A_p$. 
It follows that the weighted  space $W^{1,p}_{\rm loc}(\mathcal C;|y|^a|z|^{-b}dz)\subset L^1_{\rm loc}(\mathcal C)$ is well defined (see for instance
 \cite[Subsection 1.9]{HKM}, where the notation $H^{1,p}_{\rm loc}(\Omega;w(z)dz)$ is used). 
In addition, thanks to \cite[Theorems 3.51, 3.66]{HKM} we have that any nonnegative and nontrivial function $u\in W^{1,p}_{\rm loc}(\mathcal C;|y|^a|z|^{-b}dz)$
satisfying
\begin{equation}
\label{eq:super}
-{\rm div}\big(|y|^a|z|^{-b}\,|\nabla u|^{p-2}\nabla u\big)\ge 0\qquad \text{in $\mathcal C$}
\end{equation}
is lower semicontinuous and positive in $\mathcal C$.

The best constant $\mathcal M({\omega_\mathcal C})$ is attained by a function 
$\Phi\in W^{1,p}_0({\omega_\mathcal C};|\Pi\sigma|^ad\sigma)$, due to the compactness of 
embedding 
$W^{1,p}_0(\omega_\mathcal C;|\Pi\sigma|^ad\sigma)\hookrightarrow L^p(\omega_\mathcal C;|\Pi\sigma|^ad\sigma)$ given by Corollary \ref{C:HS_sphere}. 
By  a standard 
argument, $\Phi$ can not change sign in ${\omega_\mathcal C}$. Thus, we can assume that $\Phi$ is nonnegative.

We use spherical coordinates to define 
\begin{equation}
\label{eq:U}
U(r\sigma)=r^{-\rmH}\Phi(\sigma)
\end{equation}
 for $r>0$ and $\sigma\in\omega_\mathcal C$. 
Since 
$$
|\nabla U|^2=(\partial_rU)^2+r^{-2}|\nabla_{\!\sigma}U|^2=r^{-2(\rmH+1)}\big(\rmH^2 \Phi^2+|\nabla_{\!\sigma} \Phi|^2\big),$$
we have that  $U\in W^{1,p}_{\rm loc}(\mathcal C;|y|^a|z|^{-b}dz)$. Moreover, for any fixed 
$v\in C^\infty_c({\mathcal C})$ it holds that
$$
\nabla U\cdot\nabla v=\partial_rU\partial_r{v}+r^{-2}\nabla_{\!\sigma}U\cdot\nabla_{\!\sigma}{v}= 
r^{-(\rmH+2)}\big(-\rmH\Phi r\partial_r {v}+\nabla_{\!\sigma}\Phi\cdot\nabla_{\!\sigma}{v}\big)
$$
Notice that $d+a-b= p(\rmH+1)$, see (\ref{eq:def_H}). We have
\begin{multline*}
\int\limits_{{\mathcal C}}|y|^a|z|^{-b}|\nabla U|^{p-2}\nabla U\cdot\nabla {v}\,dz\\
=\int\limits_{{\omega_\mathcal C}}d\sigma \int\limits_0^\infty 
 |\Pi\sigma|^a \bigl(\rmH^2\Phi^2+|\nabla_{\!\sigma}\Phi|^2\big)^{\frac{p-2}{2}}
(-\rmH  \Phi~\! r^{\rmH}~\!\partial_r {v}+ r^{\rmH-1}\nabla_{\!\sigma}\Phi\cdot\nabla_{\!\sigma}{v})
~\! dr\\
=\int\limits_0^\infty r^{\rmH-1} dr\int\limits_{{\omega_\mathcal C}}
 |\Pi\sigma|^a \bigl(\rmH^2\Phi^2+|\nabla_{\!\sigma}\Phi|^2\big)^{\frac{p-2}{2}}(\rmH^2\Phi {v}+\nabla_{\!\sigma}\Phi\cdot\nabla_{\!\sigma}{v})
~\!d\sigma,
\end{multline*}
where we used integration by parts and then Fubini's theorem.

Since $\Phi$ achieves $\mathcal M(\omega_\mathcal C)$, we infer that
$$
\begin{aligned}
\int\limits_{{\mathcal C}}|y|^a|z|^{-b}&|\nabla U|^{p-2}\nabla U\cdot\nabla {v}\,dz
=\mathcal M(\omega_\mathcal C)\int\limits_0^\infty r^{\rmH-1}dr\int\limits_{\omega_\mathcal C} 
|\Pi\sigma|^a \Phi^{p-1}{v}~\!d\sigma\\
&= \mathcal M(\omega_\mathcal C)
\int\limits_0^\infty r^{p\rmH-1}dr \int\limits_{{\omega_\mathcal C}} |\Pi\sigma|^a U^{p-1}{v}~\!d\sigma =
\mathcal M(\omega_\mathcal C)
\int\limits_{{\mathcal C}}|y|^a|z|^{-b-p}U^{p-1}{v}~\!dz~\!.
\end{aligned}
$$
We proved that $U$ 
is a nonnegative local solution to
\begin{equation}
\label{eq:EL}
-{\rm div}\big(|y|^a|z|^{-b}\,|\nabla u|^{p-2}\nabla u\big)=\mathcal M({\omega_\mathcal C})\ |y|^a|z|^{-b-p}\, |u|^{p-2}u
\quad\text{in $\mathcal C$,}
\end{equation}
in the sense that
\begin{equation}
\label{eq:weak}
\int\limits_{{\mathcal C}}|y|^a|z|^{-b}|\nabla U|^{p-2}\nabla U\cdot\nabla {v}\,dz=
\mathcal M(\omega_\mathcal C)
\int\limits_{{\mathcal C}}|y|^a|z|^{-b-p}U^{p-1}{v}~\!dz
\end{equation}
for any ${v}\in C^\infty_c({\mathcal C})$. 
Since $U\in W^{1,p}_{\rm loc}(\mathcal C;|y|^a|z|^{-b}dz)$,
then (\ref{eq:weak}) holds for any $v\in {\mathcal D}^{1,p}_0(\mathcal C;|y|^a|z|^{-b}dz)$ with compact support.

To prove the inequality $\mathcal M({\omega_\mathcal C})\le m(\mathcal C)$ 
fix $v\in C^\infty_c(\mathcal C)$. 
Since $U\in W^{1,p}_{\rm loc}(\mathcal C;|y|^a|z|^{-b}dz)$ solves (\ref{eq:super}), then $U$ is bounded away
from zero on the support of $v$. It follows that 
$U^{1-p}|v|^p  \in {\mathcal D}_0^{1,p}(\mathcal C;|y|^a|z|^{-b}dz)$, hence it can be used as test function in 
(\ref{eq:weak}). We infer the equality
\begin{equation}
\label{eq:aligned}
\int\limits_{{\mathcal C}}|y|^a|z|^{-b}|\nabla U|^{p-2}\nabla U\cdot\nabla 
(U^{1-p}|v|^p)\,dz=
\mathcal M(\omega_\mathcal C)
\int\limits_{{\mathcal C}}|y|^a|z|^{-b-p}|v|^p~\!dz.
\end{equation}
To handle the first integral in (\ref{eq:aligned}) we notice that 
$$
\nabla U\cdot\nabla (U^{1-p}|v|^p)=p U^{1-p}(\nabla U\cdot\nabla v) |v|^{p-2}v
-(p-1)U^{-p}|\nabla U|^2|v|^p.
$$
Thus
$$
|\nabla U|^{p-2}\nabla U\cdot\nabla (U^{1-p}|v|^p)\le p\Big(\frac{|\nabla U||v|}{U}\Big)^{p-1}|\nabla v|-(p-1)
\Big(\frac{|\nabla U||v|}{U}\Big)^{p}\le 
|\nabla v|^p
$$
by the elementary Young inequality $p|s|^{p-1}|t|\le |t|^p+(p-1)|s|^p$. 
Thus we  conclude that
$$
\int\limits_{\mathcal C}|y|^a|z|^{-b}|\nabla v|^p~\!dz\ge   \mathcal 
M(\omega_\mathcal C)
\int\limits_{{\mathcal C}}|y|^a|z|^{-b-p}~\!|v|^p~\!dz.
$$
Since $v$ was arbitrarily chosen in $C^\infty_c(\mathcal C)$, we proved that 
$\mathcal M({\omega_\mathcal C})\le m(\mathcal C)$.

\medskip

To prove the opposite inequality we adopt a standard strategy. Consider the 
sequence
$$
u_\delta(r,\sigma)=r^{-\rmH\pm \delta}\Phi(\sigma)\qquad
\text{on $\mathcal C_\pm$,}
$$
where we have set
$\mathcal C_+=\mathcal C\cap B_1$, $\mathcal C_-=\mathcal C\setminus B_1$.
Clearly, $u_\delta\in {\mathcal D}^{1,p}_0({\mathcal C};|y|^a|z|^{-b} dz)$ and moreover
\begin{equation}
\label{eq:compare}
\int\limits_{\mathcal C}|y|^a|z|^{-b-p}|u_\delta|^p~\!dz=
\int\limits_{\mathcal C_+}|y|^a|z|^{-b-p}|u_\delta|^p~\!dz+\int\limits_{\mathcal C_-}|y|^a|z|^{-b-p}|u_\delta|^p~\!dz=
\frac{2}{p\delta}\int\limits_{\omega_\mathcal C}|\Pi\sigma|^a|\Phi|^p~\!d\sigma.
\end{equation}
It is easy to see that $|\nabla u_\delta|=r^{-\rmH-1\pm\delta}\Big[(\rmH\mp\delta)^2\Phi^2+|\nabla_{\!\sigma}\Phi|^2\Big]^\frac12$ on $\mathcal C_\pm$,
from which one easily infer
$$
\begin{aligned}
\int\limits_{\mathcal C}|y|^a|z|^{-b}|\nabla u_\delta|^p~\!dz&=\frac{1}{p\delta}
\int\limits_{\omega_\mathcal C}|\Pi\sigma|^a\big\{\big[(\rmH-\delta)^2\Phi^2+|\nabla_{\!\sigma}\Phi|^2\big]^\frac{p}{2}+
\big[(\rmH+\delta)^2\Phi^2+|\nabla_{\!\sigma}\Phi|^2\big]^\frac{p}{2}\big\}~\!d\sigma\\
&=\frac{2}{p\delta}\Big(\int\limits_{\omega_\mathcal C}|\Pi\sigma|^a\big[\rmH^2\Phi^2+|\nabla_{\!\sigma}\Phi|^2\big]^\frac{p}{2}~\!d\sigma+O(\delta^2)\Big)
\end{aligned}
$$
as $\delta\to 0$. Since $\Phi$ achieves $\mathcal M({\omega_\mathcal C})$, we infer that
$$
\int\limits_{\mathcal C}|y|^a|z|^{-b}|\nabla u_\delta|^p~\!dz=\frac{2}{p\delta}
\Big(\mathcal M({\omega_\mathcal C})~\!\int\limits_{\omega_\mathcal C}|\Pi\sigma|^a|\Phi|^p~\!d\sigma+O(\delta^2)\Big)~\!.
$$
Taking into account the definition of $m(\mathcal C)$ and
(\ref{eq:compare}), we see that 
$$
m(\mathcal C)~\le ~
\frac{\displaystyle\int_{\mathcal C}|y|^a|z|^{-b}|\nabla u_\delta|^p~\!dz}{\displaystyle\int_{\mathcal C}|\Pi\sigma|^a|u_\delta|^p~\!d\sigma}=
\mathcal M({\omega_\mathcal C})+O(\delta^2),
$$
which concludes the proof.  \QED

\paragraph{Proof of Theorem \ref{T:not_achieved}.}
As in the proof of Theorem \ref{T:cone_sphere1} we omit the indexes $p,a,b$.

If $m(\mathcal C)=0$ the result is trivial, as nonzero constant functions are not in 
${\mathcal D}^{1,p}_0(\mathcal C;|y|^a|z|^{-b}dz)$. Thus, let $m(\mathcal C)>0$.

Arguing as in the proof of Theorem \ref{T:cone_sphere1}, we can assume that (\ref{eq:stronger}) is satisfied.
Let  $U(r\sigma)=r^{-\rmH}\Phi(\sigma)$ be the function in (\ref{eq:U}).  We already proved that 
$U\in W^{1,p}_{\rm loc}(\mathcal C; |y|^a|z|^{-b}dz)$ is a lower semicontinuous and  positive local solution to 
\begin{equation}
\label{eq:ELm}
-{\rm div}\big(|y|^a|z|^{-b}\,|\nabla u|^{p-2}\nabla u\big)=m(\mathcal C)\ |y|^a|z|^{-b-p}\, |u|^{p-2}u
\quad\text{in $\mathcal C$,}
\end{equation}
compare with
(\ref{eq:EL}) and recall that $m(\mathcal C)=\mathcal M({\omega_\mathcal C})$.
From (\ref{eq:ELm}) we infer that $U$ is locally bounded outside $\Sigma_0$ by \cite{Ser}, therefore it is
of class $C^{1,\alpha}$ on $\mathcal C\setminus\Sigma_0$ by \cite{Di, To}.

By contradiction, assume that 
$u\in {\mathcal D}^{1,p}_0(\mathcal C;|y|^a|z|^{-b}dz)$ achieves $m(\mathcal C)$. Up to a change of sign,
$u$ is nonnegative and is a weak solution to (\ref{eq:ELm}). Thus, just as $U$, the function 
$u$ is positive on $\mathcal C$ and of class $C^{1,\alpha}$ outside the singular set $\Sigma_0$.

Take a domain $A$, compactly contained in $\mathcal C\setminus\Sigma_0$.
The main step in the proof consists in showing that the open sets
$$A_1=\big\{z\in A~|~~ \nabla u\cdot\nabla U \lneq |\nabla u||\nabla U| \big\}~,\quad
A_2=\Big\{z\in A~|~~
\frac{|\nabla u|}{u}\neq \frac{|\nabla U|}{U}~\Big\}
$$
are empty. This easily imply that there exists a constant $\lambda>0$ such that $u=\lambda U$ on $A$. 
Since $A$ was arbitrarily chosen,  the constant $\lambda$ does not depend on $A$. Thus $u$ is proportional to $U$ on $\mathcal C\setminus\Sigma_0$,
which is impossible as $U\notin {\mathcal D}^{1,p}_0(\mathcal C;|y|^a|z|^{-b}dz)$.

To show that the sets $A_1, A_2$ are empty we refine the calculations in the proof of 
Theorem \ref{T:cone_sphere1}. 
Take a sequence of nonnegative functions $u_h\in C^\infty_c(\mathcal C)$ such that $u_h\to u$ in $\mathcal D^{1,p}_0(\mathcal C; |y|^a|z|^{-b}dz)$.
We have
\begin{multline}
\label{eq:A1}
m(\mathcal C)
\int\limits_{{\mathcal C}}|y|^a|z|^{-b-p}u^p~\!dz=
m(\mathcal C)
\int\limits_{{\mathcal C}}|y|^a|z|^{-b-p}~\!u_h^p~\!dz+o_h(1)\\
\le 
\int\limits_{{\mathcal C}}|y|^a|z|^{-b}|\nabla U|^{p-2}
\Big(p~\!\frac{\nabla U\cdot\nabla u_h}{U^{p-1}} u_h^{p-1}
-(p-1)\frac{|\nabla U|^2u_h^p}{U^p}
\Big)dz+o_h(1)~\!.
\end{multline}
Now, as $h\to \infty$ we have
\begin{multline*}
\nu_1:=p\int\limits_{A_1}|y|^a|z|^{-b}|\nabla U|^{p-2}
\Big(\frac{|\nabla U|~\!|\nabla u|}{U^{p-1}}-\frac{\nabla U\cdot\nabla u}{U^{p-1}}\Big) ~\!u^{p-1}~\!dz\\=
p\int\limits_{A_1}|y|^a|z|^{-b}|\nabla U|^{p-2}
\frac{|\nabla U|~\!|\nabla u_h|-\nabla U\cdot\nabla u_h}{U^{p-1}} ~\!u_h^{p-1}~\!dz+o_h(1)~\!.
\end{multline*}
Notice that the nonnegative constant $\nu_1$ is positive if and only if $A_1$ has positive measure. 
From (\ref{eq:A1}) it follows that 
\begin{multline*}
m(\mathcal C)
\int\limits_{{\mathcal C}}|y|^a|z|^{-b-p}u^p~\!dz\\
\le 
\int\limits_{{\mathcal C}}|y|^a|z|^{-b}
\bigg[p\Big(\frac{|\nabla U|u_h}{U}\Big)^{p-1}|\nabla u_h|
-(p-1)
\Big(\frac{|\nabla U|u_h}{U}\Big)^{p}\bigg]dz - \nu_1+o_h(1).
\end{multline*}

Next, as $h\to \infty$ we have 
\begin{multline*}
\nu_2:=
\int\limits_{A_2}|y|^a|z|^{-b}
\Big(|\nabla u|^p+(p-1)\Big(\frac{|\nabla U|u}{U}\Big)^p-
p\Big(\frac{|\nabla U|u}{U}\Big)^{p-1}|\nabla u|\Big)dz\\
=\int\limits_{A_2}|y|^a|z|^{-b}
\Big(|\nabla u_h|^p+(p-1)\Big(\frac{|\nabla U|u_h}{U}\Big)^p-
p\Big(\frac{|\nabla U|u_h}{U}\Big)^{p-1}|\nabla u_h|\Big)dz+o_h(1)~\!.
\end{multline*}
Thanks to Young's inequality, we have that $\nu_2>0$ if and only if $A_2$ has positive measure.

We proved that
$$
m(\mathcal C)
\int\limits_{{\mathcal C}}|y|^a|z|^{-b-p}~\!u^p~\!dz
\le 
\int\limits_{{\mathcal C}}|y|^a|z|^{-b}|\nabla u_h|^p~\!dz - (\nu_1+\nu_2)+o_h(1),
$$
and letting $h\to \infty$ we infer that 
$$
\begin{aligned}
m(\mathcal C)
\int\limits_{{\mathcal C}}|y|^a|z|^{-b-p}~\!u^p~\!dz
&\le 
\int\limits_{{\mathcal C}}|y|^a|z|^{-b}|\nabla u|^p~\!dz - (\nu_1+\nu_2)\\
&=m(\mathcal C)
\int\limits_{{\mathcal C}}|y|^a|z|^{-b-p}~\!u^p~\!dz - (\nu_1+\nu_2)~\!,
\end{aligned}
$$
as $u$ achieves $m(\mathcal C)$. This implies that $\nu_1=\nu_2=0$, which is equivalent to say that 
the sets $A_1, A_2$ are empty, as claimed.
\QED

\paragraph{Proof of Theorem \ref{T:strict1}.}
We can find an open geodesic ball $\mathcal B_\theta\subset \S^{d-1}$ such that $\mathcal B_\theta\cap \partial\omega_\mathcal C\neq \emptyset$
and $\overline{\mathcal B_\theta}\cap \Sigma_0=\emptyset$. 

The restriction operator
$W^{1,p}_0(\omega_\mathcal C;|\Pi\sigma|^ad\sigma)\to W^{1,p}(\mathcal B_\theta)$ is continuous, as
$|\Pi\sigma|$ is smooth and bounded away from $0$ on $\mathcal B_\theta$.
On the other hand, any function $\f\in W^{1,p}_0(\omega_\mathcal C;|\Pi\sigma|^ad\sigma)$ vanishes on the open set 
$\mathcal B_\theta\setminus\overline{\omega_\mathcal C}$, which is not empty. We infer that 
$W^{1,p}_0(\omega_\mathcal C;|\Pi\sigma|^ad\sigma)$ does not contain constant functions,
which is enough, thanks to Theorem \ref{T:cone_sphere1}.
\QED

\paragraph{Proof of Theorem \ref{T:strict3}.} The last claim in Theorem \ref{T:cone_sphere1} gives  
$m_{2,a,b}(\R^d\setminus\Sigma_0)=\rmH^2_{2,a,b}$ if $k+a\ge 2$. Thus we can assume  $k+a<2$. In view of Theorem \ref{T:cone_sphere1}, we have
\begin{equation}
\label{eq:last}
m_{2,a,b}(\R^d\setminus\Sigma_0)=\rmH^2_{2,a,b}~+ \min_{\f\in  W^{1,2}_0(\S^{d-1}\setminus\Sigma_0;|\Pi\sigma|^ad\sigma)\atop \f\neq 0}
\frac{\displaystyle\int_{\S^{d-1}}|\Pi\sigma|^a|\nabla_{\!\sigma}\f|^2~\!d\sigma}{\displaystyle\int_{\S^{d-1}} |\Pi\sigma|^a|\f|^2~\!d\sigma}~\!.
\end{equation}
Since the embedding $W^{1,2}_0(\S^{d-1}\setminus\Sigma_0;|\Pi\sigma|^ad\sigma)\to 
L^2(\S^{d-1};|\Pi\sigma|^ad\sigma)$ 
is compact, see Corollary \ref{C:HS_sphere}, the minimum in (\ref{eq:last}) is the first eigenvalue $\lambda_1$ of the problem
\begin{equation}
\label{eq:cos}
\begin{cases}
-\text{div}_\sigma(|\Pi\sigma|^a\nabla_{\!\sigma} \f)=\lambda |\Pi\sigma|^a \f&\text{in $\S^{d-1}\setminus\Sigma_0$}\\
\f\in W^{1,2}_0(\S^{d-1}\setminus\Sigma_0;|\Pi\sigma|^ad\sigma).
\end{cases}
\end{equation}
By known facts, $\lambda_1$ is simple and the corresponding eigenfunction $\f_1$ is the only nonnegative one.
By direct computations based on the remarks above, one can  check that 
$$
\f_1(\sigma)=|\Pi\sigma|^{2-(k+a)}~,\qquad \lambda_1=(d-k)(2-(k+a)).
$$
In fact, if we write $|y|=r\cos(\theta)$, 
$|x|=r\sin(\theta)$, $0<\theta<\frac{\pi}2$, then for $\f$ depending only on 
$\theta$, problem (\ref{eq:cos}) is rewritten as follows:
\begin{equation}
\label{eq:cos2}
-(\cos^{k+a-1}(\theta)\sin^{d-k-1}(\theta) \f_\theta)_\theta =\lambda 
\cos^a(\theta)\sin^{d-k-1}(\theta) \f.
\end{equation}
So, it is evident that $\f_1=\cos^{2-a-k}(\theta)$ satisfies $(\ref{eq:cos2})$ with 
$\lambda=\lambda_1$, which concludes the proof. 
\QED

\end{document}